\documentclass[reqno]{amsart}

\usepackage{amsmath}
\usepackage{amsthm}
\usepackage{amsfonts}
\usepackage{amssymb}
\usepackage{latexsym}
\usepackage{amscd}
\usepackage{stmaryrd}
\usepackage{color}
\usepackage[all]{xypic}
\usepackage{epsfig}
\usepackage{graphics}
\usepackage{ifthen}
\usepackage{varioref}
\usepackage{rotating}

\theoremstyle{plain}
\newtheorem{thm}{Theorem}[section]
\newtheorem{cor}[thm]{Corollary}
\newtheorem{lem}[thm]{Lemma}

\newtheorem{conj}[thm]{Conjecture}
\newtheorem{exam}[thm]{Example}

\numberwithin{equation}{section}

\definecolor{darkgreen}{rgb}{0.0625,0.64,0.0625}
\usepackage[
            pdfstartview=FitH,
            CJKbookmarks=true,
            bookmarksnumbered=true,
            bookmarksopen=true,
            colorlinks,
            pdfborder=001,
            linkcolor=blue,
            anchorcolor=green,
            citecolor=red
            ]{hyperref}

\newfont{\scyr}{wncyr10 scaled 550}

\def\proof{\noindent {\bf Proof.\;}}

\allowdisplaybreaks

\begin{document}

\title{Weighted sum formulas of multiple zeta values with even arguments}

\date{\today\thanks{The first author is supported by the National Natural Science Foundation of
China (Grant No. 11471245) and Shanghai Natural Science Foundation (grant no. 14ZR1443500).} }

\author{Zhonghua Li \quad and \quad Chen Qin}

\address{School of Mathematical Sciences, Tongji University, No. 1239 Siping Road,
Shanghai 200092, China}

\email{zhonghua\_li@tongji.edu.cn}

\address{School of Mathematical Sciences, Tongji University, No. 1239 Siping Road,
Shanghai 200092, China}

\email{2014chen\_qin@tongji.edu.cn}



\begin{abstract}
We prove a weighted sum formula of the zeta values at even arguments, and a weighted sum formula of the multiple zeta values with even arguments and its zeta-star analogue. The weight coefficients are given by (symmetric) polynomials of the arguments. These weighted sum formulas for the zeta values and for the multiple zeta values were conjectured by L. Guo, P. Lei and J. Zhao.
\end{abstract}

\maketitle




\section{Introduction}\label{Sec:Intro}

For a positive integer $n$ and a sequence $\mathbf{k}=(k_1,\ldots,k_n)$ of positive integers with $k_1>1$, the multiple zeta value $\zeta(\mathbf{k})$ and the multiple zeta-star value $\zeta^{\star}(\mathbf{k})$ are defined by the following infinite series
$$\zeta(\mathbf{k})=\zeta(k_1,\ldots,k_n)=\sum\limits_{m_1>\cdots>m_n\geqslant 1}\frac{1}{m_1^{k_1}\cdots m_n^{k_n}}$$
and
$$\zeta^{\star}(\mathbf{k})=\zeta^{\star}(k_1,\ldots,k_n)=\sum\limits_{m_1\geqslant \cdots\geqslant m_n\geqslant 1}\frac{1}{m_1^{k_1}\cdots m_n^{k_n}},$$
respectively. The number $n$ is called the depth. In depth one case, both $\zeta(\mathbf{k})$ and $\zeta^{\star}(\mathbf{k})$ are special values of the Riemann zeta function at positive integer arguments.

The study of these values may be traced back to L. Euler. Among other things, L. Euler found the following sum formula
$$\sum\limits_{i=2}^{k-1}\zeta^{\star}(i,k-i)=(k-1)\zeta(k),\quad k\geqslant 3,$$
or equivalently,
$$\sum\limits_{i=2}^{k-1}\zeta(i,k-i)=\zeta(k),\quad k\geqslant 3.$$
There are many generalizations and variations of the sum formula, among which we mention some weighted sum formulas at even arguments. In \cite{Gangl-Kaneko-Zagier}, the following formula
$$\sum\limits_{i=1}^{k-1}\zeta(2i,2k-2i)=\frac{3}{4}\zeta(2k)$$
was proved by using the regularized double shuffle relations of the double zeta values. M. E. Hoffman considered the sum
$$\sum\limits_{k_1+\cdots+k_n=k\atop k_j\geqslant 1}\zeta(2k_1,\ldots,2k_n)$$
in \cite{Hoffman2017}, and we showed in \cite{Li-Qin} that the formulas given by M. E. Hoffman in \cite{Hoffman2017} are consequences of the regularized double shuffle relations of the multiple zeta values. Later in \cite{Guo-Lei-Zhao}, new families of weighted sum formulas of the forms
$$
\sum\limits_{k_1+k_2=k\atop k_j\geqslant 1}F(k_1,k_2)\zeta(2k_1)\zeta(2k_2),\quad \sum\limits_{k_1+k_2+k_3=k\atop k_j\geqslant 1}G(k_1,k_2,k_3)\zeta(2k_1)\zeta(2k_2)\zeta(2k_3)
$$
and
$$\sum\limits_{k_1+k_2=k\atop k_j\geqslant 1}F(k_1,k_2)\zeta(2k_1,2k_2),\quad \sum\limits_{k_1+k_2+k_3=k\atop k_j\geqslant 1}G(k_1,k_2,k_3)\zeta(2k_1,2k_2,2k_3)$$
were given, where $F(x,y)$ and $G(x,y,z)$ are (symmetric) polynomials with rational coefficients. And in the end of \cite{Guo-Lei-Zhao}, L. Guo, P. Lei and J. Zhao proposed the following general conjecture.

\begin{conj}[{\cite[Conjecture 4.7]{Guo-Lei-Zhao}}]\label{Conj:WeightedSum}
Let $F(x_1,\ldots,x_n)\in\mathbb{Q}[x_1,\ldots,x_n]$ be a symmetric polynomial of degree $r$. Set $d=\deg_{x_1}F(x_1,\ldots,x_n)$. Then for every positive integer $k\geqslant n$ we have
\begin{align}
\sum\limits_{k_1+\cdots+k_n=k\atop k_j\geqslant 1}F(k_1,\ldots,k_n)\zeta(2k_1)\cdots\zeta(2k_n)=\sum\limits_{l=0}^Te_{F,l}(k)\zeta(2l)\zeta(2k-2l),
\label{Eq:WeightedSum-Zeta-Conj}\\
\sum\limits_{k_1+\cdots+k_n=k\atop k_j\geqslant 1}F(k_1,\ldots,k_n)\zeta(2k_1,\ldots,2k_n)=\sum\limits_{l=0}^Tc_{F,l}(k)\zeta(2l)\zeta(2k-2l),
\label{Eq:WeightedSum-MZV-Conj}
\end{align}
where $T=\max\{[(r+n-2)/2],[(n-1)/2]\}$, $e_{F,l}(x),c_{F,l}(x)\in\mathbb{Q}[x]$ depend only on $l$ and $F$, $\deg e_{F,l}(x)\leqslant r-1$ and $\deg c_{F,l}(x)\leqslant d$.
\end{conj}

Here as usual, for a real number $x$, we denote by $[x]$ the greatest integer that not exceeding $x$.

The purpose of this paper is to give a proof of Conjecture \ref{Conj:WeightedSum}. In fact, we also prove the zeta-star analogue of \eqref{Eq:WeightedSum-MZV-Conj}. To prove \eqref{Eq:WeightedSum-Zeta-Conj}, as in \cite{Guo-Lei-Zhao}, we first establish a weighted sum formula of the Bernoulli numbers. In \cite{Guo-Lei-Zhao}, L. Guo, P. Lei and J. Zhao used certain zeta functions to study the Bernoulli numbers. Here we just use the generating function of the Bernoulli numbers. Hence our method seems more elementary. After getting the weighted sum formula of the Bernoulli numbers, we obtain the weighted sum formula \eqref{Eq:WeightedSum-Zeta-Conj} by using Euler's evaluation formula of the zeta values at even arguments. Finally, applying the symmetric sum formulas of M. E. Hoffman \cite{Hoffman1992}, we obtain the weighted sum formula \eqref{Eq:WeightedSum-MZV-Conj} and its zeta-star analogue from the formula \eqref{Eq:WeightedSum-Zeta-Conj}.

The paper is organized as follows. In Section \ref{Sec:WeightSum-Bernoulli}, we deal with the weighted sum of the Bernoulli numbers. In Section \ref{Sec:WeightSum-Zeta}, we prove the weighted sum formula \eqref{Eq:WeightedSum-Zeta-Conj}. And in Section \ref{Sec:WeightSum-MZV}, we prove the weighted sum formula \eqref{Eq:WeightedSum-MZV-Conj} and its zeta-star analogue. Finally, in Section \ref{Sec:RegDouble-WeightSum}, we show that the weighted sum formulas obtained in this paper can be deduced from the regularized double shuffle relations of the multiple zeta values.


\section{A weighted sum formula of Bernoulli numbers}\label{Sec:WeightSum-Bernoulli}

The Bernoulli numbers $\{B_i\}$ are defined by
$$\sum\limits_{i=0}^{\infty} \frac{B_i}{i!}t^i=\frac{t}{e^t-1}.$$
It is known that $B_0=1$, $B_1=-\frac{1}{2}$ and $B_i=0$ for odd $i\geqslant 3$. We set
\begin{align*}
&f(t)=\frac{t}{e^t-1}-1+\frac{1}{2}t=\sum\limits_{i=1}^\infty\frac{B_{2i}}{(2i)!}t^{2i},\\
&g(t)=\frac{t}{e^t-1}+\frac{1}{2}t=\sum\limits_{i=0}^\infty\frac{B_{2i}}{(2i)!}t^{2i}.
\end{align*}
We compute the derivatives of the even function $f(t)$. Let $D=t\frac{d}{dt}$ and $h(t)=\frac{t}{e^t-1}$. Then using the formula
$$h'(t)=\frac{1-t}{e^t-1}-\frac{t}{(e^t-1)^2},$$
we find that for any nonnegative integer $m$,
\begin{align}
D^mf(t)=\sum\limits_{i=0}^{m+1}f_{mi}(t)h(t)^i.
\label{Eq:Diff-f}
\end{align}
Here $f_{mi}(t)$ are polynomials determined by $f_{00}(t)=\frac{1}{2}t-1$, $f_{01}(t)=1$ and the recursive formulas
\begin{align}
\begin{cases}
f_{m0}(t)=tf_{m-1,0}'(t) & \text{for\;} m\geqslant 1,\\
f_{m,m+1}(t)=-mf_{m-1,m}(t) & \text{for\;} m\geqslant 1,\\
f_{mi}(t)=tf'_{m-1,i}(t)+i(1-t)f_{m-1,i}(t)-(i-1)f_{m-1,i-1}(t) & \text{for\;} 1\leqslant i\leqslant m.
\end{cases}
\label{Eq:Recursive-fmi}
\end{align}
In particular, for any integers $m,i$ with $1\leqslant i\leqslant m+1$, we have $f_{mi}(t)\in\mathbb{Z}[t]$.
From \eqref{Eq:Recursive-fmi}, it is easy to see that for any nonnegative integer $m$, we have
$$f_{m0}(t)=\frac{1}{2}t-\delta_{m,0},\quad f_{m,m+1}(t)=(-1)^mm!.$$

\begin{lem}
For any integers $m,i$ with $1\leqslant i\leqslant m+1$, we have $\deg f_{mi}(t)=m+1-i$, and the leading coefficient $c_{mi}$ of $f_{mi}(t)$ satisfies the condition $(-1)^mc_{mi}>0$.
\end{lem}

\proof We use induction on $m$. Assume that $m\geqslant 1$. The result for $i=m+1$ follows from $f_{m,m+1}(t)=(-1)^mm!$. Now assume the integer $i$ satisfies the condition $1\leqslant i\leqslant m$, and
$$f_{m-1,i}(t)=c_{m-1,i}t^{m-i}+\text{lower degree terms}$$
with $(-1)^{m-1}c_{m-1,i}>0$.
Let $c_{m0}=\frac{1}{2}$. Then we have
$$f_{mi}(t)=(-ic_{m-1,i}-(i-1)c_{m-1,i-1})t^{m+1-i}+\text{lower degree terms}.$$
As
\begin{align*}
&(-1)^m(-ic_{m-1,i}-(i-1)c_{m-1,i-1})\\
=&i(-1)^{m-1}c_{m-1,i}+(i-1)(-1)^{m-1}c_{m-1,i-1}>0,
\end{align*}
we get the result.
\qed

Therefore we have
$$f_{m0}(t)=c_{m0}t-\delta_{m,0}$$
with $c_{m0}=\frac{1}{2}$, and for any integers $m,i$ with the condition $1\leqslant i\leqslant m+1$, we have
$$f_{mi}(t)=c_{mi}t^{m+1-i}+\text{lower degree terms},$$
with the recursive formula
$$c_{mi}=-ic_{m-1,i}-(i-1)c_{m-1,i-1},\quad (1\leqslant i\leqslant m)$$
and $c_{m,m+1}=(-1)^mm!$.

\begin{cor}
For any nonnegative integer $m$, we have $c_{m1}=(-1)^m$.
\end{cor}

For later use, we need the following lemma.

\begin{lem}
For any nonnegative integer $m$, we have
\begin{align}
\sum\limits_{i=1}^{m+1}(-1)^{i-1}f_{mi}(t)t^{i-1}=1.
\label{Eq:Sum-fmi}
\end{align}
In particular, we have
\begin{align}
\sum\limits_{i=1}^{m+1}(-1)^{i-1}c_{mi}=\delta_{m,0}.
\label{Eq:Sum-cmi}
\end{align}
\end{lem}

\proof
We proceed by induction on $m$ to prove \eqref{Eq:Sum-fmi}. The case of $m=0$ follows from the fact $f_{01}(t)=1$. Now assume that $m\geqslant 1$, using the recursive formula \eqref{Eq:Recursive-fmi}, we have
\begin{align*}
&\sum\limits_{i=1}^{m+1}(-1)^{i-1}f_{mi}(t)t^{i-1}=\sum\limits_{i=1}^m(-1)^{i-1}f_{m-1,i}'(t)t^{i}+\sum\limits_{i=1}^{m}(-1)^{i-1}if_{m-1,i}(t)t^{i-1}\\
&\qquad+\sum\limits_{i=1}^m(-1)^iif_{m-1,i}(t)t^i+\sum\limits_{i=1}^m(-1)^{i}(i-1)f_{m-1,i-1}(t)t^{i-1}+m!t^m\\
=&\sum\limits_{i=1}^m(-1)^{i-1}(f_{m-1,i}(t)t^{i})'+(-1)^mmf_{m-1,m}(t)t^m+m!t^m\\
=&\frac{d}{dt}\sum\limits_{i=1}^m(-1)^{i-1}f_{m-1,i}(t)t^{i}.
\end{align*}
Then we get \eqref{Eq:Sum-fmi} from the induction assumption. Finally, comparing the coefficients of $t^m$ of both sides of \eqref{Eq:Sum-fmi}, we get \eqref{Eq:Sum-cmi}.
\qed

Now we want to express $h(t)^i$ by $D^mg(t)$. For this purpose, we use matrix computations. For any nonnegative integer $m$, let $A_m(t)$ be a $(m+1)\times (m+1)$ matrix defined by
$$A_m(t)=\begin{pmatrix}
f_{01}(t) & &&\\
f_{11}(t) & f_{12}(t) &&\\
\vdots & \vdots & \ddots &\\
f_{m1}(t) & f_{m2}(t) & \cdots & f_{m,m+1}(t)
\end{pmatrix}.$$
Note that for $m\geqslant 1$, we have
$$A_m(t)=\begin{pmatrix}
A_{m-1}(t) & 0\\
\alpha_m(t) & (-1)^mm!
\end{pmatrix}$$
with $\alpha_m(t)=(f_{m1}(t),\ldots,f_{mm}(t))$. From linear algebra, we know that the matrix $\begin{pmatrix}
A & 0\\
C & B
\end{pmatrix}$ is invertible with
$$\begin{pmatrix}
A & 0\\
C & B
\end{pmatrix}^{-1}=\begin{pmatrix}
A^{-1} & 0\\
-B^{-1}CA^{-1} & B^{-1}
\end{pmatrix},$$
provided that $A$ and $B$ are invertible square matrices. Therefore by induction on $m$, we find that for all nonnegative integer $m$, the matrices $A_m(t)$ are invertible, and the inverses satisfy the recursive formula
\begin{align}
A_m(t)^{-1}=\begin{pmatrix}
A_{m-1}(t)^{-1} & 0\\
(-1)^{m+1}\frac{1}{m!}\alpha_m(t)A_{m-1}(t)^{-1} & (-1)^m\frac{1}{m!}
\end{pmatrix},\quad (m\geqslant 1).
\label{Eq:Recursive-AmInverse}
\end{align}

For any nonnegative integer $m$, set
$$A_m(t)^{-1}=\begin{pmatrix}
g_{01}(t) & &&\\
g_{11}(t) & g_{12}(t) &&\\
\vdots & \vdots & \ddots &\\
g_{m1}(t) & g_{m2}(t) & \cdots & g_{m,m+1}(t)
\end{pmatrix}.$$

\begin{lem}
Let $m$ and $i$ be integers.
\begin{itemize}
  \item [(1)] For any $m\geqslant 0$, we have $g_{m,m+1}(t)=(-1)^m\frac{1}{m!}$;
  \item [(2)] If $1\leqslant i\leqslant m$, we have the recursive formula
  \begin{align}
  g_{mi}(t)=(-1)^{m+1}\frac{1}{m!}\sum\limits_{j=i}^{m}f_{mj}(t)g_{j-1,i}(t);
  \label{Eq:Recursive-gmi}
  \end{align}
  \item [(3)] If $1\leqslant i\leqslant m+1$, we have $g_{mi}(t)\in\mathbb{Q}[t]$ with $\deg g_{mi}(t)\leqslant m+1-i$;
  \item [(4)] For $1\leqslant i\leqslant m+1$, set
  $$g_{mi}(t)=d_{mi}t^{m+1-i}+\text{lower degree terms}.$$
  Then we have $d_{m,m+1}=(-1)^m\frac{1}{m!}$ and
   \begin{align}
  d_{mi}=(-1)^{m+1}\frac{1}{m!}\sum\limits_{j=i}^{m}c_{mj}d_{j-1,i}
  \label{Eq:Recursive-dmi}
  \end{align}
  for $1\leqslant i\leqslant m$.
\end{itemize}
\end{lem}

\proof The assertions in items (1) and (2) follow from \eqref{Eq:Recursive-AmInverse}. To prove the item (3), we proceed by induction on $m$. For the case of $m=0$, we get the result from $g_{01}(t)=1$. Assume that $m\geqslant 1$, then $g_{m,m+1}(t)=(-1)^m\frac{1}{m!}\in\mathbb{Q}[t]$ with degree zero. For $1\leqslant i\leqslant j\leqslant m$, using the induction assumption, we may set
$$g_{j-1,i}(t)=d_{j-1,i}t^{j-i}+\text{lower degree terms}\in\mathbb{Q}[t].$$
Since
$$f_{mj}(t)=c_{mj}t^{m+1-j}+\text{lower degree terms}\in\mathbb{Z}[t],$$
we get
$$f_{mj}(t)g_{j-1,i}(t)=c_{mj}d_{j-1,i}t^{m+1-i}+\text{lower degree terms}\in\mathbb{Q}[t].$$
Using \eqref{Eq:Recursive-gmi}, we finally get
$$g_{mi}(t)=\left((-1)^{m+1}\frac{1}{m!}\sum\limits_{j=i}^{m}c_{mj}d_{j-1,i}\right)t^{m+1-i}+\text{lower degree terms}\in\mathbb{Q}[t].$$
The item (4) follows from the above proof.
\qed

\begin{cor}
For any nonnegative integer $m$, we have $d_{m1}=(-1)^m$.
\end{cor}

\proof
We use induction on $m$. If $m\geqslant 1$, using \eqref{Eq:Recursive-dmi} and the induction assumption, we get
$$d_{m1}=(-1)^{m+1}\frac{1}{m!}\sum\limits_{j=1}^m(-1)^{j-1}c_{mj}.$$
By \eqref{Eq:Sum-cmi}, we have
$$d_{m1}=(-1)^{m+1}\frac{1}{m!}(\delta_{m,0}-(-1)^{m}c_{m,m+1}),$$
which implies the result.
\qed

To get $h(t)^i$, we rewrite \eqref{Eq:Diff-f} as
\begin{align}
\begin{pmatrix}
g(t)\\
Dg(t)\\
\vdots\\
D^mg(t)
\end{pmatrix}-\frac{1}{2}t\begin{pmatrix}
1\\
1\\
\vdots\\
1
\end{pmatrix}=A_m(t)\begin{pmatrix}
h(t)\\
h(t)^2\\
\vdots\\
h(t)^{m+1}
\end{pmatrix},
\label{Eq:G-h-Matrix}
\end{align}
and rewrite \eqref{Eq:Sum-fmi} as
$$A_{m}(t)\begin{pmatrix}
1\\
-t\\
t^2\\
\vdots\\
(-1)^mt^m
\end{pmatrix}=\begin{pmatrix}
1\\
1\\
\vdots\\
1
\end{pmatrix}.$$
Therefore we find
$$
\begin{pmatrix}
h(t)\\
h(t)^2\\
\vdots\\
h(t)^{m+1}
\end{pmatrix}=A_m(t)^{-1}\begin{pmatrix}
g(t)\\
Dg(t)\\
\vdots\\
D^mg(t)
\end{pmatrix}-\frac{1}{2}t\begin{pmatrix}
1\\
-t\\
t^2\\
\vdots\\
(-1)^mt^m
\end{pmatrix}.
$$
Then for any positive integer $i$, we get
\begin{align}
h(t)^i=\sum\limits_{j=1}^ig_{i-1,j}(t)D^{j-1}g(t)+\frac{1}{2}(-1)^it^i.
\label{Eq:h-i}
\end{align}

For the later use, we prepare a lemma.

\begin{lem}\label{Lem:G-linearInd}
For a nonnegative integer $m$, the functions $1,g(t),D g(t),\ldots,D^m g(t)$ are linearly independent over the rational function field $\mathbb{Q}(t)$. \end{lem}

\proof
Let $p(t),p_0(t),p_1(t),\ldots,p_m(t)\in\mathbb{Q}(t)$ satisfy
$$p(t)+p_0(t)g(t)+p_1(t)D g(t)+\cdots+p_m(t)D^mg(t)=0.$$
Using \eqref{Eq:G-h-Matrix}, we get
$$p(t)+\frac{1}{2}t\sum\limits_{j=0}^mp_j(t)+(p_0(t),\ldots,p_m(t))A_m(t)\begin{pmatrix}
h(t)\\
\vdots\\
h(t)^{m+1}
\end{pmatrix}=0.$$
Since $e^t$ is transcendental over $\mathbb{Q}(t)$, we know $e^t-1$, and then $h(t)$ is transcendental over $\mathbb{Q}(t)$. Hence we have
$$p(t)+\frac{1}{2}t\sum\limits_{j=0}^mp_j(t)=0,\qquad (p_0(t),\ldots,p_m(t))A_m(t)=0,$$
which implies that all  $p_j(t)$ and $p(t)$ are zero functions as the matrix $A_m(t)$ is invertible.
\qed

From now on let $n$ be a fixed positive integer, and $m_1,\ldots,m_n$ be fixed nonnegative integers. We want to compute $D^{m_1}f(t)\cdots D^{m_n}f(t)$. On the one hand, using \eqref{Eq:Diff-f}, we have
$$D^{m_1}f(t)\cdots D^{m_n}f(t)=\sum\limits_{i=0}^{m_1+\cdots+m_n+n}f_i(t)h(t)^i,$$
with
$$f_i(t)=\sum\limits_{i_1+\cdots+i_n=i\atop 0\leqslant i_j\leqslant m_j+1}f_{m_1i_1}(t)\cdots f_{m_ni_n}(t).$$

\begin{lem}
We have
$$f_0(t)=\prod\limits_{j=1}^n\left(\frac{1}{2}t-\delta_{m_j,0}\right),$$
and $\deg f_i(t)\leqslant m_1+\cdots+m_n+n-i$ for any nonnegative integer $i$.
\end{lem}

\proof
For integers $i_1,\ldots,i_n$ with the conditions $i_1+\cdots+i_n=i$ and $0\leqslant i_j\leqslant m_j+1$, we have
$$\deg(f_{m_1i_1}(t)\cdots f_{m_ni_n}(t))\leqslant \sum\limits_{j=1}^n(m_j+1-i_j)=m_1+\cdots+m_n+n-i,$$
which deduces that $\deg f_i(t)\leqslant m_1+\cdots+m_n+n-i$.
\qed

Then using \eqref{Eq:h-i}, we get
\begin{align*}
&D^{m_1}f(t)\cdots D^{m_n}f(t)\\
=&\sum\limits_{i=1}^{m_1+\cdots+m_n+n}f_i(t)\left(\sum\limits_{j=1}^ig_{i-1,j}(t)D^{j-1}g(t)+\frac{1}{2}(-1)^it^i\right)+f_0(t)\\
=&\sum\limits_{j=1}^{m_1+\cdots+m_n+n}F_j(t)D^{j-1}g(t)+F_0(t)
\end{align*}
with
$$F_0(t)=f_0(t)+\frac{1}{2}\sum\limits_{i=1}^{m_1+\cdots+m_n+n}(-1)^if_i(t)t^i$$
and
$$F_j(t)=\sum\limits_{i=j}^{m_1+\cdots+m_n+n}f_i(t)g_{i-1,j}(t),\quad (1\leqslant j\leqslant m_1+\cdots+m_n+n).$$

\begin{lem}
Let $j$ be a nonnegative integer with $j\leqslant m_1+\cdots+m_n+n$. Then
\begin{itemize}
  \item [(1)] the function $F_j(t)$ is even;
  \item [(2)] we have
  $$F_0(t)=\frac{1}{2}\prod\limits_{j=1}^n\left(\frac{1}{2}t-\delta_{m_j,0}\right)+\frac{1}{2}(-1)^n
\prod\limits_{j=1}^n\left(\frac{1}{2}t+\delta_{m_j,0}\right).$$
In particular, $\deg F_0(t)\leqslant n$;
  \item [(3)] if $j>0$, we have $\deg F_j(t)\leqslant m_1+\cdots+m_n+n-j$.
Moreover, we have $\deg F_1(t)\leqslant m_1+\cdots+m_n+n-2$ provided that $n$ is even or $m_1,\ldots,m_n$ are not all zero.
\end{itemize}
\end{lem}

\proof
Since $D^{m}f(t)$ and $D^m g(t)$ are even, we have
$$\sum\limits_{j=1}^{m_1+\cdots+m_n+n}F_j(t)D^{j-1}g(t)+F_0(t)=\sum\limits_{j=1}^{m_1+\cdots+m_n+n}F_j(-t)D^{j-1}g(t)+F_0(-t).$$
Then by Lemma \ref{Lem:G-linearInd}, we know all $F_j(t)$ are even functions.

By the definition of $f_i(t)$, we have
$$\sum\limits_{i=0}^{m_1+\cdots+m_n+n}(-1)^if_i(t)t^i=\prod\limits_{j=1}^n\sum\limits_{i_j=0}^{m_j+1}(-1)^{i_j}f_{m_ji_j}(t)t^{i_j}.$$
Using \eqref{Eq:Sum-fmi}, we find
$$\sum\limits_{i=0}^{m_1+\cdots+m_n+n}(-1)^if_i(t)t^i=\prod\limits_{j=1}^n(f_{m_j0}(t)-t).$$
Then we get (2) from the fact that $f_{m0}(t)=\frac{1}{2}t-\delta_{m,0}$ and the expression of $f_0(t)$.

Since
$$\deg f_i(t)g_{i-1,j}(t)\leqslant (m_1+\cdots+m_n+n-i)+(i-j)=m_1+\cdots+m_n+n-j,$$
we get $\deg F_j(t)\leqslant m_1+\cdots+m_n+n-j$.

If we set
$$\widetilde{c}_{mi}=\begin{cases}
\frac{1}{2}\delta_{m,0} & \text{if\;} i=0,\\
c_{mi} & \text{if\;} i\neq 0,
\end{cases}$$
then the coefficient of $t^{m+1-i}$ in $f_{mi}(t)$ is $\widetilde{c}_{mi}$ for any integers $m,i$ with the condition $0\leqslant i\leqslant m+1$.
Since
$$F_1(t)=\sum\limits_{i=1}^{m_1+\cdots+m_n+n}\sum\limits_{i_1+\cdots+i_n=i\atop 0\leqslant i_j\leqslant m_j+1}f_{m_1i_1}(t)\cdots f_{m_ni_n}(t)g_{i-1,1}(t),$$
and $d_{i-1,1}=(-1)^{i-1}$, we find the coefficient of $t^{m_1+\cdots+m_n+n-1}$ in $F_1(t)$ is
\begin{align*}
&\sum\limits_{i=1}^{m_1+\cdots+m_n+n}\sum\limits_{i_1+\cdots+i_n=i\atop 0\leqslant i_j\leqslant m_j+1}(-1)^{i-1}\widetilde{c}_{m_1i_1}\cdots \widetilde{c}_{m_ni_n}\\
=&\widetilde{c}_{m_10}\cdots \widetilde{c}_{m_n0}-\prod\limits_{j=1}^n\sum\limits_{i_j=0}^{m_j+1}(-1)^{i_j}\widetilde{c}_{m_ji_j},
\end{align*}
which is
$$\widetilde{c}_{m_10}\cdots \widetilde{c}_{m_n0}-\prod\limits_{j=1}^n(\widetilde{c}_{m_j0}-\delta_{m_j,0})$$
by \eqref{Eq:Sum-cmi}. Then the coefficient of $t^{m_1+\cdots+m_n+n-1}$ in $F_1(t)$ is
$$\left(\frac{1}{2}\right)^n(1-(-1)^n)\delta_{m_1,0}\cdots\delta_{m_n,0},$$
which is zero if $n$ is even or at least one $m_i$ is not zero.
\qed

Now for a positive integer $j$ with $j\leqslant m_1+\cdots+m_n+n$, let $a_{jl}\in\mathbb{Q}$ be the coefficient of $t^{2l}$ in the polynomial $F_j(t)$. Then we have
\begin{align}
F_j(t)=\sum\limits_{l=0}^{\left[\frac{m_1+\cdots+m_n+n-j}{2}\right]}a_{jl}t^{2l}.
\label{Eq:Fj}
\end{align}
If $n$ is even or $m_1,\ldots,m_n$ are not all zero, we have
$$F_1(t)=\sum\limits_{l=0}^{\left[\frac{m_1+\cdots+m_n+n-2}{2}\right]}a_{1l}t^{2l}.$$
Hence we have
$$D^{m_1}f(t)\cdots D^{m_n}f(t)=\sum\limits_{j=1}^{m_1+\cdots+m_n+n}\sum\limits_{l=0}^{\left[\frac{m_1+\cdots+m_n+n-j}{2}\right]}a_{jl}t^{2l}D^{j-1}g(t)+F_0(t).$$
Changing the order of the summation, we have
$$D^{m_1}f(t)\cdots D^{m_n}f(t)=\sum\limits_{l=0}^{T}\sum\limits_{j=1}^{m_1+\cdots+m_n+n-2l}a_{jl}t^{2l}D^{j-1}g(t)+F_0(t),$$
where
$$T=\begin{cases}
\left[\frac{n-1}{2}\right] & \text{if\;} m_1=\cdots=m_n=0,\\
\left[\frac{m_1+\cdots+m_n+n-2}{2}\right] & \text{otherwise}.
\end{cases}$$
Since
$$D^{j-1}g(t)=\sum\limits_{i=0}^\infty(2i)^{j-1}\frac{B_{2i}}{(2i)!}t^{2i},$$
we get
\begin{align*}
&D^{m_1}f(t)\cdots D^{m_n}f(t)\\
=&\sum\limits_{k=0}^\infty\sum\limits_{l=0}^{\min\{T,k\}}\left(\sum\limits_{j=1}^{m_1+\cdots+m_n+n-2l}a_{jl}(2k-2l)^{j-1}\right)\frac{B_{2k-2l}}{(2k-2l)!}t^{2k}+F_0(t).
\end{align*}
Then the coefficient of $t^{2k}$ in $D^{m_1}f(t)\cdots D^{m_n}f(t)$ is
\begin{align}
\sum\limits_{l=0}^{\min\{T,k\}}\left(\sum\limits_{j=1}^{m_1+\cdots+m_n+n-2l}2^{j-1}a_{jl}(k-l)^{j-1}\right)\frac{B_{2k-2l}}{(2k-2l)!},
\label{Eq:Coeff-Right}
\end{align}
provided that $k\geqslant n$.

On the other hand, since
$$D^mf(t)=\sum\limits_{i=1}^{\infty}(2i)^m\frac{B_{2i}}{(2i)!}t^{2i},$$
we find the coefficient of $t^{2k}$ in $D^{m_1}f(t)\cdots D^{m_n}f(t)$ is
\begin{align}
\sum\limits_{k_1+\cdots+k_n=k\atop k_j\geqslant 1}(2k_1)^{m_1}\cdots(2k_n)^{m_n}\frac{B_{2k_1}\cdots B_{2k_n}}{(2k_1)!\cdots(2k_n)!}.
\label{Eq:Coeff-Left}
\end{align}

Finally, comparing \eqref{Eq:Coeff-Right} with \eqref{Eq:Coeff-Left}, we get a weighted sum formula of the Bernoulli numbers.

\begin{thm}\label{Thm:WeightedSum-Bernoulli}
Let $n,k$ be positive integers with $k\geqslant n$. Then for any nonnegative integers $m_1,\ldots,m_n$, we have
\begin{align}
&\sum\limits_{k_1+\cdots+k_n=k\atop k_j\geqslant 1}k_1^{m_1}\cdots k_n^{m_n}\frac{B_{2k_1}\cdots B_{2k_n}}{(2k_1)!\cdots(2k_n)!}\nonumber\\
=&\sum\limits_{l=0}^{\min\{T,k\}}\left(\sum\limits_{j=1}^{m_1+\cdots+m_n+n-2l}\frac{a_{jl}}{2^{m_1+\cdots+m_n-j+1}}(k-l)^{j-1}\right)\frac{B_{2k-2l}}{(2k-2l)!},
\label{Eq:WeightedSum-Bernoulli}
\end{align}
where $T=\max\{[(m_1+\cdots+m_n+n-2)/2],[(n-1)/2]\}$ and $a_{jl}$ are determined by \eqref{Eq:Fj}.
\end{thm}

Note that in \cite[Theorem 1]{Petojevic-Srivastava}, A. Petojevi\'{c} and H. M. Srivastava had considered the case of $m_1=\cdots=m_n=0$. See also \cite[Theorems 1 and 2]{Dilcher}.

In the end of this section, we list some explicit examples of $n=4$. Note that some examples of $n=2$ and $n=3$ were given in \cite{Guo-Lei-Zhao}.

\begin{exam}\label{Exe:Bernoulli}
Let $k$ be a positive integer with $k\geqslant 4$. Set $\sum=\sum\limits_{k_1+k_2+k_3+k_4=k\atop k_j\geqslant 1}$. We have
\begin{align*}
&\sum\frac{B_{2k_1}B_{2k_2}B_{2k_3}B_{2k_4}}{(2k_1)!(2k_2)!(2k_3)!(2k_4)!}=-\frac{(k+1)(2k+1)(2k+3)}{3}\frac{B_{2k}}{(2k)!}
-\frac{2k}{3}\frac{B_{2k-2}}{(2k-2)!},\\
&\sum k_1^2\frac{B_{2k_1}B_{2k_2}B_{2k_3}B_{2k_4}}{(2k_1)!(2k_2)!(2k_3)!(2k_4)!}=-\frac{k(k+1)(2k+1)(2k+3)(4k+3)}{120}\frac{B_{2k}}{(2k)!}\\
&\qquad\qquad-\frac{k(4k^2-6k+3)}{24}\frac{B_{2k-2}}{(2k-2)!}-\frac{2k-5}{160}\frac{B_{2k-4}}{(2k-4)!},\\
&\sum k_1^3\frac{B_{2k_1}B_{2k_2}B_{2k_3}B_{2k_4}}{(2k_1)!(2k_2)!(2k_3)!(2k_4)!}=-\frac{k(k+1)(2k+1)(2k+3)(4k^2+6k+1)}{240}
\frac{B_{2k}}{(2k)!}\\
&\qquad-\frac{k(12k^3-12k^2-11k+9)}{96}\frac{B_{2k-2}}{(2k-2)!}-\frac{(2k-5)(13k-9)}{960}\frac{B_{2k-4}}{(2k-4)!}.
\end{align*}
Set $B_{k_1,k_2,k_3,k_4}=\frac{B_{2k_1}B_{2k_2}B_{2k_3}B_{2k_4}}{(2k_1)!(2k_2)!(2k_3)!(2k_4)!}$. Using the formulas
\begin{align*}
&\sum k_1 B_{k_1,k_2,k_3,k_4}
=\frac{k}{4}\sum B_{k_1,k_2,k_3,k_4},\\
&\sum k_1k_2 B_{k_1,k_2,k_3,k_4}
=\frac{k^2}{12}\sum B_{k_1,k_2,k_3,k_4}-\frac{1}{3}\sum k_1^2 B_{k_1,k_2,k_3,k_4},\\
&\sum k_1^2k_2 B_{k_1,k_2,k_3,k_4}
=\frac{k}{3}\sum k_1^2B_{k_1,k_2,k_3,k_4}-\frac{1}{3}\sum k_1^3 B_{k_1,k_2,k_3,k_4},\\
&\sum k_1k_2k_3 B_{k_1,k_2,k_3,k_4}
=\frac{k^3}{24}\sum B_{k_1,k_2,k_3,k_4}-\frac{k}{2}\sum k_1^2 B_{k_1,k_2,k_3,k_4}\\
&\qquad\qquad\qquad\qquad+\frac{1}{3}\sum k_1^3 B_{k_1,k_2,k_3,k_4},
\end{align*}
one can work out other weighted sum formulas of the Bernoulli numbers with the condition $m_1+m_2+m_3+m_4\leqslant 3$.
\end{exam}


\section{Weighted sum formulas of zeta values at even arguments}\label{Sec:WeightSum-Zeta}

Euler's formula claims that for any positive integer $k$,
\begin{align}
\zeta(2k)=(-1)^{k+1}\frac{B_{2k}}{2(2k)!}(2\pi)^{2k}.
\label{Eq:Euler-Formula}
\end{align}
Then from Theorem \ref{Thm:WeightedSum-Bernoulli}, we get the following weighted sum formula for zeta values at even arguments.

\begin{thm}\label{Thm:WeightedSum-Zeta}
Let $n,k$ be positive integers with $k\geqslant n$. Then for any nonnegative integers $m_1,\ldots,m_n$, we have
\begin{align}
&\sum\limits_{k_1+\cdots+k_n=k\atop k_j\geqslant 1}k_1^{m_1}\cdots k_n^{m_n}\zeta(2k_1)\cdots\zeta(2k_n)=(-1)^n\sum\limits_{l=0}^{\min\{T,k\}}\frac{(2l)!}{B_{2l}}\nonumber\\
&\times\left(\sum\limits_{j=1}^{m_1+\cdots+m_n+n-2l}\frac{a_{jl}}{2^{m_1+\cdots+m_n+n-j-1}}(k-l)^{j-1}\right)\zeta(2l)\zeta(2k-2l),
\label{Eq:WeightedSum-Zeta}
\end{align}
where $T=\max\{[(m_1+\cdots+m_n+n-2)/2],[(n-1)/2]\}$ and $a_{jl}$ are determined by \eqref{Eq:Fj}.
\end{thm}

Finally, we obtain the weighted sum formula \eqref{Eq:WeightedSum-Zeta-Conj}.

\begin{thm}\label{Thm:WeightedSum-Zeta-General}
Let $n,k$ be positive integers with $k\geqslant n$. Let $F(x_1,\ldots,x_n)\in\mathbb{Q}[x_1,\ldots,x_n]$ be a polynomial of degree $r$. Then we have
$$\sum\limits_{k_1+\cdots+k_n=k\atop k_j\geqslant 1}F(k_1,\ldots,k_n)\zeta(2k_1)\cdots \zeta(2k_n)=\sum\limits_{l=0}^{\min\{T,k\}}e_{F,l}(k)\zeta(2l)\zeta(2k-2l),$$
where $T=\max\{[(r+n-2)/2],[(n-1)/2]\}$, $e_{F,l}(x)\in\mathbb{Q}[x]$ depends only on $l$ and $F$, and $\deg e_{F,l}(x)\leqslant r+n-2l-1$.
\end{thm}

Note that the polynomial $F(x_1,\ldots,x_n)$ in Theorem \ref{Thm:WeightedSum-Zeta-General} need not be symmetric, and the upper bound for the degree of the polynomial $e_{F,l}(x)$ is different from that in Conjecture \ref{Conj:WeightedSum}.  In Conjecture \ref{Conj:WeightedSum}, the upper bound for $\deg e_{F,l}(x)$ is $r-1$, which should be a typo. See the examples below.

\begin{exam}
Set $\sum=\sum\limits_{k_1+k_2+k_3+k_4=k\atop k_j\geqslant 1}$. For a positive integer $k$ with $k\geqslant 4$, we have
\begin{align*}
&\sum \zeta(2k_1)\zeta(2k_2)\zeta(2k_3)\zeta(2k_4)=\frac{(k+1)(2k+1)(2k+3)}{24}\zeta(2k)
-2k\zeta(2)\zeta(2k-2),\\
&\sum k_1^2\zeta(2k_1)\zeta(2k_2)\zeta(2k_3)\zeta(2k_4)=\frac{k(k+1)(2k+1)(2k+3)(4k+3)}{960}\zeta(2k)\\
&\qquad\qquad-\frac{k(4k^2-6k+3)}{8}\zeta(2)\zeta(2k-2)+\frac{9(2k-5)}{8}\zeta(4)\zeta(2k-4),\\
&\sum k_1^3\zeta(2k_1)\zeta(2k_2)\zeta(2k_3)\zeta(2k_4)=\frac{k(k+1)(2k+1)(2k+3)(4k^2+6k+1)}{1920}
\zeta(2k)\\
&\quad-\frac{k(12k^3-12k^2-11k+9)}{32}\zeta(2)\zeta(2k-2)+\frac{3(2k-5)(13k-9)}{16}\zeta(4)\zeta(2k-4),
\end{align*}
which can deduce all other weighted sums \eqref{Eq:WeightedSum-Zeta} under the conditions $n=4$ and $m_1+m_2+m_3+m_4\leqslant 3$ as explained in Example \ref{Exe:Bernoulli}.
\end{exam}


\section{Weighted sum formulas of multiple zeta values with even arguments}\label{Sec:WeightSum-MZV}

To treat the weighted sum of the multiple zeta values with even arguments and its zeta-star analogue, we recall the symmetric sum formulas of M. E. Hoffman \cite[Theorems 2.1 and 2.2]{Hoffman1992}. For a partition $\Pi=\{P_1,P_2,\ldots,P_i\}$ of the set $\{1,2,\ldots,n\}$, let $l_j=\sharp P_j$ and
$$c(\Pi)=\prod\limits_{j=1}^i (l_j-1)!,\quad \tilde{c}(\Pi)=(-1)^{n-i}c(\Pi).$$
We also denote by $\mathcal{P}_n$ the set of all partitions of the set $\{1,2,\ldots,n\}$. Then the symmetric sum formulas are
\begin{align}
\sum\limits_{\sigma\in S_n}\zeta(k_{\sigma(1)},\ldots,k_{\sigma(n)})=\sum\limits_{\Pi\in\mathcal{P}_n}\tilde{c}(\Pi)\zeta(\mathbf{k},\Pi)
\label{Eq:SymSum-MZV}
\end{align}
and
\begin{align}
\sum\limits_{\sigma\in S_n}\zeta^{\star}(k_{\sigma(1)},\ldots,k_{\sigma(n)})=\sum\limits_{\Pi\in\mathcal{P}_n}c(\Pi)\zeta(\mathbf{k},\Pi),
\label{Eq:SymSum-MZSV}
\end{align}
where $\mathbf{k}=(k_1,\ldots,k_n)$ is a sequence of positive integers with all $k_i>1$, $S_n$ is the symmetric group of degree $n$ and for a partition $\Pi=\{P_1,\ldots,P_i\}\in\mathcal{P}_n$,
$$\zeta(\mathbf{k},\Pi)=\prod\limits_{j=1}^i\zeta\left(\sum\limits_{l\in P_j}k_l\right).$$

Now let $\mathbf{k}=(2k_1,\ldots,2k_n)$ with all $k_i$ positive integers. Using \eqref{Eq:SymSum-MZV} and \eqref{Eq:SymSum-MZSV}, we have
\begin{align}
&\sum\limits_{\sigma\in S_n}\zeta(2k_{\sigma(1)},\ldots,2k_{\sigma(n)})\nonumber\\
=&\sum\limits_{i=1}^n(-1)^{n-i}\sum\limits_{l_1+\cdots+l_i=n\atop l_1\geqslant \cdots\geqslant l_i\geqslant 1}\prod\limits_{j=1}^i(l_j-1)!\sum\limits_{\Pi=\{P_1,\ldots,P_i\}\in\mathcal{P}_n\atop \sharp{P_j}=l_j}\zeta(\mathbf{k},\Pi)
\label{Eq:SymSum-2-MZV}
\end{align}
and
\begin{align}
&\sum\limits_{\sigma\in S_n}\zeta^{\star}(2k_{\sigma(1)},\ldots,2k_{\sigma(n)})\nonumber\\
=&\sum\limits_{i=1}^n\sum\limits_{l_1+\cdots+l_i=n\atop l_1\geqslant \cdots\geqslant l_i\geqslant 1 }\prod\limits_{j=1}^i(l_j-1)!\sum\limits_{\Pi=\{P_1,\ldots,P_i\}\in\mathcal{P}_n\atop \sharp{P_j}=l_j}\zeta(\mathbf{k},\Pi).
\label{Eq:SymSum-2-MZSV}
\end{align}

From now on, let $k,n$ be fixed positive integers with $k\geqslant n$, and let $F(x_1,\ldots,x_n)$ be a fixed symmetric polynomial with rational coefficients. It is easy to see that
\begin{align*}
&\sum\limits_{k_1+\cdots+k_n=k\atop k_j\geqslant 1}F(k_1,\ldots,k_n)\sum\limits_{\sigma\in S_n}\zeta(2k_{\sigma(1)},\ldots,2k_{\sigma(n)})\\
=&n!\sum\limits_{k_1+\cdots+k_n=k\atop k_j\geqslant 1}F(k_1,\ldots,k_n)\zeta(2k_1,\ldots,2k_n)
\end{align*}
and
\begin{align*}
&\sum\limits_{k_1+\cdots+k_n=k\atop k_j\geqslant 1}F(k_1,\ldots,k_n)\sum\limits_{\sigma\in S_n}\zeta^{\star}(2k_{\sigma(1)},\ldots,2k_{\sigma(n)})\\
=&n!\sum\limits_{k_1+\cdots+k_n=k\atop k_j\geqslant 1}F(k_1,\ldots,k_n)\zeta^{\star}(2k_1,\ldots,2k_n).
\end{align*}
On the other hand, for a partition $\Pi=\{P_1,\ldots,P_i\}\in\mathcal{P}_n$ with $\sharp P_j=l_j$, we have
\begin{align}
&\sum\limits_{k_1+\cdots+k_n=k\atop k_j\geqslant 1}F(k_1,\ldots,k_n)\zeta(\mathbf{k},\Pi)\nonumber\\
=&\sum\limits_{t_1+\cdots+t_i=k\atop t_j\geqslant 1}\sum\limits_{{{k_1+\cdots+k_{l_1}=t_1\atop\vdots}\atop k_{l_1+\cdots+l_{i-1}+1}+\cdots+k_n=t_i}\atop k_j\geqslant 1}F(k_1,\ldots,k_n)\zeta(2t_1)\cdots\zeta(2t_i).
\label{Eq:F-times-zeta}
\end{align}
To treat the inner sum about $F(k_1,\ldots,k_n)$ in the right-hand side of \eqref{Eq:F-times-zeta}, we need the following lemmas.

\begin{lem}\label{Lem:PowerSum-2}
For any positive integer $k$ and any nonnegative integers $p_1,p_2$, we have
\begin{align}
&\sum\limits_{i=1}^{k-1} i^{p_1}(k-i)^{p_2}\nonumber\\
=&\sum\limits_{i,j\geqslant 0\atop p_1\leqslant i+j\leqslant p_1+p_2}(-1)^{j+p_1}\binom{i+j}{i}\binom{p_2}{i+j-p_1}\frac{B_i}{j+1}(k-1)^{j+1}k^{p_1+p_2-i-j}.
\label{Eq:PowerSum-2}
\end{align}
In particular, the right-hand side of \eqref{Eq:PowerSum-2} is a polynomial of $k$ with rational coefficients of degree $p_1+p_2+1$.
\end{lem}

\proof
Let $S_{p_1,p_2}(k)=\sum\limits_{i=1}^{k-1} i^{p_1}(k-i)^{p_2}$ and let
$$G_k(t_1,t_2)=\sum\limits_{p_1,p_2\geqslant 0}S_{p_1,p_2}(k)\frac{t_1^{p_1}t_2^{p_2}}{p_1!p_2!}$$
be the generating function. We have
$$G_k(t_1,t_2)=\sum\limits_{i=1}^{k-1}e^{it_1+(k-i)t_2}=\frac{(1-e^{(k-1)(t_1-t_2)})e^{kt_2}}{e^{t_2-t_1}-1}.$$
Using the definition of the Bernoulli numbers, we get
\begin{align*}
G_k(t_1,t_2)=&\sum\limits_{i\geqslant 0,j\geqslant 1,l\geqslant 0}(-1)^{i}\frac{B_i}{i!j!l!}(k-1)^jk^l(t_1-t_2)^{i+j-1}t_2^l\\
=&\sum\limits_{i,j,l\geqslant 0}(-1)^{i}\frac{B_i}{i!(j+1)!l!}(k-1)^{j+1}k^l(t_1-t_2)^{i+j}t_2^l.
\end{align*}
Finally, we obtain the expansion
$$G_k(t_1,t_2)=\sum\limits_{i,j,l\geqslant 0\atop 0\leqslant m\leqslant i+j}(-1)^{j+m}\binom{i+j}{m}\frac{B_i}{i!(j+1)!l!}(k-1)^{j+1}k^lt_1^mt_2^{i+j+l-m}.$$
Comparing the coefficient of $\frac{t_1^{p_1}t_2^{p_2}}{p_1!p_2!}$, we get \eqref{Eq:PowerSum-2}.

Then as a polynomial of $k$, the degree of the right-hand side of \eqref{Eq:PowerSum-2} is less than or equal to $p_1+p_2+1$, and the coefficient of $k^{p_1+p_2+1}$ is
$$\sum\limits_{j=p_1}^{p_1+p_2}(-1)^{j+p_1}\binom{p_2}{j-p_1}\frac{1}{j+1},$$
which is
$$\sum\limits_{j=0}^{p_2}(-1)^{j}\binom{p_2}{j}\frac{1}{j+p_1+1}.$$
Since
$$x^{p_1}(1-x)^{p_2}=\sum\limits_{j=0}^{p_2}(-1)^{j}\binom{p_2}{j}x^{j+p_1},$$
we find the coefficient of $k^{p_1+p_2+1}$ is
$$\int_0^1x^{p_1}(1-x)^{p_2}dx=B(p_1+1,p_2+1)=\frac{p_1!p_2!}{(p_1+p_2+1)!},$$
which is nonzero.
\qed

More generally, we have

\begin{lem}\label{Lem:PowerSum}
Let $k$ and $n$ be integers with $k\geqslant n\geqslant 1$, and let $p_1,\ldots,p_n$ be nonnegative integers.  Then there exists a polynomial $f(x)\in\mathbb{Q}[x]$ of degree $p_1+\cdots+p_n+n-1$, such that
$$\sum\limits_{k_1+\cdots+k_n=k\atop k_j\geqslant 1}k_1^{p_1}\cdots k_n^{p_n}=f(k).$$
\end{lem}

\proof
We proceed by induction on $n$. If $n=1$, we may take $f(x)=x^{p_1}$. For $n>1$, since
$$\sum\limits_{k_1+\cdots+k_n=k\atop k_j\geqslant 1}k_1^{p_1}\cdots k_n^{p_n}=\sum\limits_{k_1+k_2=k\atop k_j\geqslant 1}\left(\sum\limits_{l_1+\cdots+l_{n-1}=k_1\atop l_j\geqslant 1}l_1^{p_1}\cdots l_{n-1}^{p_{n-1}}\right)k_2^{p_n},$$
using the induction assumption we have
$$\sum\limits_{k_1+\cdots+k_n=k\atop k_j\geqslant 1}k_1^{p_1}\cdots k_n^{p_n}=\sum\limits_{k_1+k_2=k\atop k_j\geqslant 1}g(k_1)k_2^{p_n},$$
where $g(x)\in\mathbb{Q}[x]$ is of degree $p_1+\cdots+p_{n-1}+n-2$. Then the result follows from Lemma \ref{Lem:PowerSum-2}.
\qed

Now we return to the computation of the right-hand side of \eqref{Eq:F-times-zeta}. Using Lemma \ref{Lem:PowerSum}, there exists a polynomial $f_{t_1,\ldots,t_i}(x_1,\ldots,x_i)\in\mathbb{Q}[x_1,\ldots,x_i]$ of degree $\deg F+n-i$, such that
\begin{align*}
&\sum\limits_{k_1+\cdots+k_n=k\atop k_j\geqslant 1}F(k_1,\ldots,k_n)\zeta(\mathbf{k},\Pi)\\
=&\sum\limits_{t_1+\cdots+t_i=k\atop t_j\geqslant 1}f_{t_1,\ldots,t_i}(t_1,\ldots,t_i)\zeta(2t_1)\cdots\zeta(2t_i).
\end{align*}
Therefore  we get
\begin{align*}
&\sum\limits_{k_1+\cdots+k_n=k\atop k_j\geqslant 1}F(k_1,\ldots,k_n)\zeta(2k_1,\ldots,2k_n)
=\frac{1}{n!}\sum\limits_{i=1}^n(-1)^{n-i}\sum\limits_{l_1+\cdots+l_i=n\atop l_1\geqslant\cdots \geqslant l_i\geqslant 1}\\
&\times\prod\limits_{j=1}^i(l_j-1)!n(l_1,\ldots,l_i)\sum\limits_{t_1+\cdots+t_i=k\atop t_j\geqslant 1}f_{t_1,\ldots,t_i}(t_1,\ldots,t_i)\zeta(2t_1)\cdots\zeta(2t_i)
\end{align*}
and
\begin{align*}
&\sum\limits_{k_1+\cdots+k_n=k\atop k_j\geqslant 1}F(k_1,\ldots,k_n)\zeta^{\star}(2k_1,\ldots,2k_n)
=\frac{1}{n!}\sum\limits_{i=1}^n\sum\limits_{l_1+\cdots+l_i=n\atop l_1\geqslant\cdots\geqslant l_i\geqslant 1}\\
&\times\prod\limits_{j=1}^i(l_j-1)!n(l_1,\ldots,l_i)\sum\limits_{t_1+\cdots+t_i=k\atop t_j\geqslant 1}f_{t_1,\ldots,t_i}(t_1,\ldots,t_i)\zeta(2t_1)\cdots\zeta(2t_i),
\end{align*}
where
$$n(l_1,\ldots,l_i)=\frac{n!}{\prod\limits_{j=1}^il_j!\prod\limits_{j=1}^n\sharp\{m\mid 1\leqslant m\leqslant i,k_m=j\}!}$$
is the number of partitions $\Pi=\{P_1,\ldots,P_i\}\in\mathcal{P}_n$ with the conditions $\sharp P_j=l_j$ for $j=1,2,\ldots,i$.

Applying Theorem \ref{Thm:WeightedSum-Zeta-General}, we then prove the weighted sum formula \eqref{Eq:WeightedSum-MZV-Conj} and its zeta-star analogue.

\begin{thm}\label{Thm:WeightedSum-MZV}
Let $n,k$ be positive integers with $k\geqslant n$. Let $F(x_1,\ldots,x_n)\in\mathbb{Q}[x_1,\ldots,x_n]$ be a symmetric polynomial of degree $r$. Then we have
$$\sum\limits_{k_1+\cdots+k_n=k\atop k_j\geqslant 1}F(k_1,\ldots,k_n)\zeta(2k_1,\ldots,2k_n)=\sum\limits_{l=0}^{\min\{T,k\}}c_{F,l}(k)\zeta(2l)\zeta(2k-2l)$$
and
$$\sum\limits_{k_1+\cdots+k_n=k\atop k_j\geqslant 1}F(k_1,\ldots,k_n)\zeta^{\star}(2k_1,\ldots,2k_n)=\sum\limits_{l=0}^{\min\{T,k\}}c_{F,l}^{\star}(k)\zeta(2l)\zeta(2k-2l),$$
where $T=\max\{[(r+n-2)/2],[(n-1)/2]\}$, $c_{F,l}(x),c_{F,l}^{\star}(x)\in\mathbb{Q}[x]$ depend only on $l$ and $F$, and $\deg c_{F,l}(x),\deg c_{F,l}^{\star}(x)\leqslant r+n-2l-1$.
\end{thm}

Note that in Theorem \ref{Thm:WeightedSum-MZV}, the upper bound for the polynomial $c_{F,l}(x)$ is different from that in Conjecture \ref{Conj:WeightedSum}. In Conjecture \ref{Conj:WeightedSum}, the upper bound for $\deg c_{F,l}(x)$ is $\deg_{x_1}F(x_1,\ldots,x_n)$. It seems that one may obtain this upper bound but need more efforts.

\begin{exam}
After getting the weighted sum formulas \eqref{Eq:WeightedSum-Zeta} with $n=2$ and $n=3$, we can obtain the weighted sum formulas of the multiple zeta values (resp. the multiple zeta-star values) of depth four. Here are some examples. For multiple zeta values, we have
\begin{align*}
&\sum\zeta(2k_1,2k_2,2k_3,2k_4)=\frac{35}{64}\zeta(2k)-\frac{5}{16}\zeta(2)\zeta(2k-2),\\
&\sum\left(k_1^2+k_2^2+k_3^2+k_4^2\right)\zeta(2k_1,2k_2,2k_3,2k_4)=\frac{7k(10k-3)}{128}\zeta(2k)\\
&\qquad-\frac{10k^2+9k-30}{32}\zeta(2)\zeta(2k-2)+\frac{3(2k-5)}{16}\zeta(4)\zeta(2k-4),\\
&\sum\left(k_1^3+k_2^3+k_3^3+k_4^3\right)\zeta(2k_1,2k_2,2k_3,2k_4)=\frac{7k(40k^2-18k+3)}{512}\zeta(2k)\\
&\qquad-\frac{40k^3+54k^2-174k+15}{128}\zeta(2)\zeta(2k-2)+\frac{3(2k-5)(3k+2)}{32}\zeta(4)\zeta(2k-4),
\end{align*}
and for multiple zeta-star values, we have
\begin{align*}
&\sum\zeta^{\star}(2k_1,2k_2,2k_3,2k_4)=\frac{(4k-5)(8k^2-20k+3)}{192}\zeta(2k)\\
&\qquad\qquad-\frac{4k-7}{16}\zeta(2)\zeta(2k-2),\\
&\sum\left(k_1^2+k_2^2+k_3^2+k_4^2\right)\zeta^{\star}(2k_1,2k_2,2k_3,2k_4)\\
=&\frac{k(128k^4-600k^3+920k^2-600k+227)}{1920}\zeta(2k)\\
&\qquad-\frac{(2k-3)(16k^2-63k+68)}{96}\zeta(2)\zeta(2k-2)-\frac{2k-5}{16}\zeta(4)\zeta(2k-4),\\
&\sum\left(k_1^3+k_2^3+k_3^3+k_4^3\right)\zeta^{\star}(2k_1,2k_2,2k_3,2k_4)\\
=&\frac{k(256k^5-1440k^4+2760k^3-2400k^2+1664k-435)}{7680}\zeta(2k)\\
&\qquad-\frac{32k^4-184k^3+318k^2-136k-51}{128}\zeta(2)\zeta(2k-2)\\
&\qquad+\frac{15(k-4)(2k-5)}{32}\zeta(4)\zeta(2k-4).
\end{align*}
Here $k$ is a positive integer with $k\geqslant 4$ and $\sum=\sum\limits_{k_1+k_2+k_3+k_4=k\atop k_j\geqslant 1}$.
\end{exam}


\section{Regularized double shuffle relations and weighted sum formulas}\label{Sec:RegDouble-WeightSum}

In this section, we briefly explain that the weighted sum formulas in Theorems \ref{Thm:WeightedSum-Zeta}, \ref{Thm:WeightedSum-Zeta-General} and \ref{Thm:WeightedSum-MZV} can be deduced from the regularized double shuffle relations of the multiple zeta values (For the details of the regularized double shuffle relations, one can refer to \cite{Ihara-Kaneko-Zagier,Racinet} or \cite{Li-Qin}).

We get Theorem \ref{Thm:WeightedSum-Zeta} and hence Theorem \ref{Thm:WeightedSum-Zeta-General} just from \eqref{Eq:WeightedSum-Bernoulli} and Euler's formula \eqref{Eq:Euler-Formula}. While \eqref{Eq:WeightedSum-Bernoulli} is an equation about the Bernoulli numbers and Euler's formula can be deduced from the regularized double shuffle relations (\cite{Li-Qin}). Hence we get Theorems \ref{Thm:WeightedSum-Zeta} and \ref{Thm:WeightedSum-Zeta-General} from the regularized double shuffle relations.

We get Theorem \ref{Thm:WeightedSum-MZV} from Theorem \ref{Thm:WeightedSum-Zeta-General} and the symmetric sum formulas. While the symmetric sum formulas are consequences of the harmonic shuffle products (\cite[Theorem 2.3]{Hoffman2015}). In fact, let $Y=\{z_k\mid k=1,2,\ldots\}$ be an alphabet with noncommutative letters and let $Y^{\ast}$ be the set of all words generated by letters in $Y$, which contains the empty word $1_Y$. Let $\mathfrak{h}^1=\mathbb{Q}\langle Y\rangle$ be the noncommutative polynomial algebra over $\mathbb{Q}$ generated by $Y$. As in \cite{Hoffman1997,Muneta}, we define two bilinear commutative products $\ast$ and $\bar{\ast}$ on $\mathfrak{h}^1$ by the rules
\begin{align*}
&1_Y\ast w=w\ast 1_Y=w,\\
&z_kw_1\ast z_lw_2=z_k(w_1\ast z_lw_2)+z_l(z_kw_1\ast w_2)+z_{k+l}(w_1\ast w_2);\\
&1_Y \,\bar{\ast}\, w=w\,\bar{\ast}\, 1_Y=w,\\
&z_kw_1\,\bar{\ast}\, z_lw_2=z_k(w_1\,\bar{\ast}\, z_lw_2)+z_l(z_kw_1\,\bar{\ast}\, w_2)-z_{k+l}(w_1\,\bar{\ast}\, w_2),
\end{align*}
where $w,w_1,w_2\in Y^{\ast}$ and $k,l$ are positive integers. Let
$$\mathfrak{h}^0=\mathbb{Q}1_Y+\sum\limits_{n,k_1,\ldots,k_n\geqslant 1\atop k_1\geqslant 2}\mathbb{Q}z_{k_1}\cdots z_{k_n}$$
be a subalgebra of $\mathfrak{h}^1$, which is also a subalgebra with respect to either the product $\ast$ or the product $\bar{\ast}$. Let $Z:\mathfrak{h}^0\rightarrow\mathbb{R}$ and $Z^{\star}:\mathfrak{h}^0\rightarrow\mathbb{R}$ be the $\mathbb{Q}$-linear maps determined by
$Z(1_Y)=Z^{\star}(1_Y)=1$ and
$$Z(z_{k_1}\cdots z_{k_n})=\zeta(k_1,\ldots,k_n),\quad Z^{\star}(z_{k_1}\cdots z_{k_n})=\zeta^{\star}(k_1,\ldots,k_n),$$
where $n,k_1,\ldots,k_n\geqslant 1$ with $k_1\geqslant 2$. It is known that both the maps $Z:(\mathfrak{h}^0,\ast)\rightarrow \mathbb{R}$ and $Z^{\star}:(\mathfrak{h}^0,\bar{\ast})\rightarrow \mathbb{R}$ are algebra homomorphisms. Hence from the following lemma, we know that the symmetric sum formulas are consequences of the harmonic shuffle products. And therefore Theorem \ref{Thm:WeightedSum-MZV} is also deduced from the regularized double shuffle relations.

\begin{lem}
Let $n$ be a positive integer and $\mathbf{k}=(k_1,\ldots,k_n)$ be a sequence of positive integers. We have
\begin{align}
\sum\limits_{\sigma\in S_n}z_{k_{\sigma(1)}}\cdots z_{k_{\sigma(n)}}=\sum\limits_{\Pi=\{P_1,\ldots,P_i\}\in\mathcal{P}_n}\tilde{c}(\Pi)z_{\mathbf{k},P_1}\ast\cdots\ast z_{\mathbf{k},P_i}
\label{Eq:SymSum-ast}
\end{align}
and
\begin{align}
\sum\limits_{\sigma\in S_n}z_{k_{\sigma(1)}}\cdots z_{k_{\sigma(n)}}=\sum\limits_{\Pi=\{P_1,\ldots,P_i\}\in\mathcal{P}_n}c(\Pi)z_{\mathbf{k},P_1}\bar{\ast}\cdots\bar{\ast} z_{\mathbf{k},P_i},
\label{Eq:SymSum-sast}
\end{align}
where $z_{\mathbf{k},P_j}=z_{\sum\limits_{l\in P_j}k_l}$.
\end{lem}

\proof
To be self contained, we give a proof here. We prove \eqref{Eq:SymSum-ast} and one can prove \eqref{Eq:SymSum-sast} similarly. We proceed by induction on $n$. The case of $n=1$ is obvious. Now assume that \eqref{Eq:SymSum-ast} is proved for $n$. Let $\mathbf{k}=(k_1,\ldots,k_n)$ and $\mathbf{k}'=(k_1,\ldots,k_{n+1})$. Since
\begin{align*}
&z_{k_{n+1}}\ast \sum\limits_{\sigma\in S_n}z_{k_{\sigma(1)}}\cdots z_{k_{\sigma(n)}}=\sum\limits_{\sigma\in S_n}\sum\limits_{j=1}^{n+1}z_{k_{\sigma(1)}}\cdots z_{k_{\sigma(j-1)}}z_{k_{n+1}}z_{k_{\sigma(j)}}\cdots z_{k_{\sigma(n)}}\\
&\qquad+\sum\limits_{\sigma\in S_n}\sum\limits_{j=1}^{n}z_{k_{\sigma(1)}}\cdots z_{k_{\sigma(j-1)}}z_{k_{\sigma(j)}+k_{n+1}}z_{k_{\sigma(j+1)}}\cdots z_{k_{\sigma(n)}}\\
=&\sum\limits_{\sigma\in S_{n+1}}z_{k_{\sigma(1)}}\cdots z_{k_{\sigma(n+1)}}+\sum\limits_{j=1}^n\sum\limits_{\sigma\in S_n}z_{k^{(j)}_{\sigma(1)}}\cdots z_{k^{(j)}_{\sigma(n)}}
\end{align*}
with
$$\mathbf{k}^{(j)}=(k_1,\ldots,k_{j-1},k_j+k_{n+1},k_{j+1},\ldots,k_n)=(k_1^{(j)},\ldots,k_n^{(j)}),$$
using the induction assumption on $\mathbf{k}$ and $\mathbf{k}^{(j)}$ with $j=1,\ldots,n$, we have
\begin{align*}
&\sum\limits_{\sigma\in S_{n+1}}z_{k_{\sigma(1)}}\cdots z_{k_{\sigma(n+1)}}=\sum\limits_{\Pi=\{P_1,\ldots,P_i\}\in\mathcal{P}_n}\tilde{c}(\Pi)z_{\mathbf{k},P_1}\ast\cdots\ast z_{\mathbf{k},P_i}\ast z_{k_{n+1}}\\
&\qquad\quad-\sum\limits_{j=1}^n\sum\limits_{\Pi=\{P_1,\ldots,P_i\}\in\mathcal{P}_n}\tilde{c}(\Pi)z_{\mathbf{k}^{(j)},P_1}\ast\cdots\ast z_{\mathbf{k}^{(j)},P_i}.
\end{align*}
Because any $\Pi\in\mathcal{P}_{n+1}$ must satisfy and can only satisfy one of the following two conditions:
\begin{itemize}
  \item [(i)] there exists one $P\in\Pi$, such that $P=\{k_{n+1}\}$;
  \item [(ii)] for any $P\in \Pi$, $P\neq \{k_{n+1}\}$,
\end{itemize}
 we see that the right-hand side of the above equation is just
$$\sum\limits_{\Pi=\{P_1,\ldots,P_i\}\in\mathcal{P}_{n+1}}\tilde{c}(\Pi)z_{\mathbf{k}',P_1}\ast\cdots\ast z_{\mathbf{k}',P_i}.$$
Hence we get \eqref{Eq:SymSum-ast}.
\qed


\end{document}